\def\ind{\par\hskip 1truecm\relax}

\let\em\it

\input diagrams

\catcode`\@=11

\newcount\secno
\newcount\subsecno
\newcount\subsubsecno
\newcount\bibno
\newdimen\tmpdimension

\newif\ifnotfound
\newif\iffound
\font\sc=cmcsc10

\def\namedef#1{\expandafter\def\csname #1\endcsname}
\def\nameuse#1{\csname #1\endcsname}
\def\typeout#1{\immediate\write16{#1}}
\long\def\ifundefined#1#2#3{\expandafter\ifx\csname#1\endcsname%
                            \relax#2\else#3\fi}
\def\hwrite#1#2{{\let\the=0\edef\next{\write#1{#2}}\next}}
\def\@ifnextchar#1#2#3{%
  \let\reserved@d=#1%
  \def\reserved@a{#2}\def\reserved@b{#3}%
  \futurelet\@let@token\@ifnch}
\def\@ifnch{%
  \ifx\@let@token\@sptoken
    \let\reserved@c\@xifnch
  \else
    \ifx\@let@token\reserved@d
      \let\reserved@c\reserved@a
    \else
      \let\reserved@c\reserved@b
    \fi
  \fi
  \reserved@c}
\let\@end\end

\def\begin#1{\ifundefined{@#1}%
                 {\def\tempa{\typeout{begin{#1} is undefined}\endgroup}}%
                 {\def\tempa{\csname @#1\endcsname}}%
             \tempa}%

\def\end#1{\ifundefined{@@#1}%
                 {\def\tempa{\typeout{end{#1} is undefined}\endgroup}}%
                 {\def\tempa{\csname @@#1\endcsname}}%
             \tempa}

\newwrite\aux
\newwrite\bbl

\def\openall{\openout\aux=\jobname.aux}

\newread\testfile
\def\lookatfile#1{\openin\testfile=\jobname.#1
    \ifeof\testfile{\immediate\openout\nameuse{#1}\jobname.#1
                    \write\nameuse{#1}{}
                    \immediate\closeout\nameuse{#1}}\fi%
    \immediate\closein\testfile}%

\def\@document{\citation{L-S:verlinde}
\citation{B-L-S:picard}
\citation{D-N:picard}
\newlabel{th:pic}{1.1}
\citation{L-S:verlinde}
\citation{L-S:verlinde}
\citation{B-L-S:picard}
\newlabel{th:local_factoriality}{1.2}
\citation{Te:BWB}
\newlabel{subsec:Lie-general-set-up}{2.1}
\newlabel{form:cent_ext}{2.1.1}
\newlabel{subsec:conf-blocks}{2.2}
\newlabel{def:conf-blocks}{2.2.1}
\newlabel{subsec:E_8}{2.3}
\citation{TUY:conf-field}
\citation{So:NB794}
\citation{So:NB794}
\citation{TUY:conf-field}
\citation{L-S:verlinde}
\citation{L-S:verlinde}
\newlabel{form:Faltings}{3.2.1}
\citation{L-S:verlinde}
\newlabel{gl(H)-ext}{3.2.2}
\newlabel{can-ext}{3.2.3}
\citation{L-S:verlinde}
\citation{L-S:verlinde}
\citation{L-S:verlinde}
\citation{L-S:verlinde}
\newlabel{tower}{3.5.1}
\citation{K-N:picard}
\citation{L-S:verlinde}
\citation{L-S:verlinde}
\citation{B-L-S:picard}
\citation{B-L-S:picard}
\citation{B-L-S:picard}
\bibdata{sorger}
\bibcite{B-L-S:picard}{B-L-S}
\bibcite{D-N:picard}{D-N}
\bibcite{K-N:picard}{K-N}
\bibcite{L-S:verlinde}{L-S}
\bibcite{So:NB794}{S}
\bibcite{Te:BWB}{T}
\bibcite{TUY:conf-field}{T-U-Y}
\bibstyle{sorger}
}

\def\@@document{\@end}

\toksdef\ta=0 \toksdef\tb=2
\long\def\leftappenditem#1\to#2{\ta={\\{#1}}\tb=\expandafter{#2}%
                                \edef#2{\the\ta\the\tb}}
\long\def\rightappenditem#1\to#2{\ta={\\{#1}}\tb=\expandafter{#2}%
                                \edef#2{\the\tb\the\ta}}

\def\lop#1\to#2{\expandafter\lopoff#1\lopoff#1#2}
\long\def\lopoff\\#1#2\lopoff#3#4{\def#4{#1}\def#3{#2}}

\def\ismember#1\of#2{\foundfalse{\let\given=#1%
    \def\\##1{\def\next{##1}%
    \ifx\next\given{\global\foundtrue}\fi}#2}}

\def\@isinlabellist#1\of#2{\notfoundtrue%
   {\def\given{#1}%
    \def\\##1{\def\next{##1}%
    \lop\next\to\za\lop\next\to\zb%
    \ifx\za\given{\zb\global\notfoundfalse}\fi}#2}%
    \ifnotfound{(?)\immediate\write16%
                 {Warning - [Page \the\pageno] {#1} No reference found}}%
                \fi}%

\def\ref#1{\ifx\@labellist\empty{(?)\immediate\write16%
                 {Warning - No references found at all.}}
               \else{\@isinlabellist{#1}\of\@labellist}\fi}

\def\newlabel#1#2{\rightappenditem{\\{#1}\\{#2}}\to\@labellist}
\def\@labellist{}

\def\label#1{\relax}

\def\cite#1{\@citenow{#1}}%

\def\bibcite#1#2{\rightappenditem{\\{#1}\\{#2}}\to\@biblist}
\def\@biblist{}

\def\@isinbiblist#1\of#2{\notfoundtrue%
   {\def\given{#1}%
    \def\\##1{\def\next{##1}%
    \lop\next\to\za\lop\next\to\zb%
    \ifx\za\given{\zb\global\notfoundfalse}\fi}#2}%
    \ifnotfound{?\immediate\write16%
                 {Warning - cite{#1} No reference found}}\fi}

\def\@citenow#1{\ifx\@biblist\empty%
                {[?]\immediate\write16{Warning - No references found}}%
                \else{[\@isinbiblist{#1}\of\@biblist]}\fi}

\def\@biblabel#1{[#1]}

\def\bibitem{\@ifnextchar[\@lbibitem\@bibitem}
\def\@lbibitem[#1]#2{\global\advance\bibno by 1%
                \item{\@biblabel{#1}}}
\def\@bibitem#1{\global\advance\bibno by 1%
                \item{\@biblabel{\number\bibno}}}

\let\newblock\relax
\def\citation#1{}
\def\bibdata#1{}
\def\bibstyle#1{}

\newbox\boxtemp
\def\@thebibliography#1%
  	{\begingroup
   	\setbox\boxtemp=\hbox{-[#1]-}
    \parindent=\wd\boxtemp
   	\global\bibno=0
   	\bigskip
   	\centerline{\bf References}
   	\medskip}
  
\def\@@thebibliography{\endgroup}

\def\bibliography#1{
}

\def\bibliographystyle#1{\relax}

\def\otstyle{article} % Par defaut
\def\@otstyle{article}
\def\section#1{\bigbreak%
               \global\def\currenvir{section}%
               \global\advance\secno by1\global\subsecno=0\global\subsubsecno=0%
               \ifx\otstyle\@otstyle\@arsection{#1}\else\@bosection{#1}\fi}
\def\@arsection#1{{\bf \ind\number\secno. {#1}}\medskip}
\def\@bosection#1{{\twelvepoint\bf \ind\number\secno. {#1}}\bigskip}

\def\subsection#1{\global\def\currenvir{subsection}%
                  \global\advance\subsecno by1\global\subsubsecno=0%
           \ifx\otstyle\@otstyle\@arsubsection{#1}\else\@bosubsection{#1}\fi}

\def\@arsubsection#1{\ind{\bf (\number\secno.\number\subsecno) }{#1}}

\def\@bosubsection#1{\medbreak\ind{\bf (\number\secno.\number\subsecno) {#1}}
                     \medskip}

\def\subsubsection#1{\medbreak\global\def\currenvir{subsection}%
                     \global\advance\subsubsecno by1%
                \ind{(\number\secno.\number\subsecno.\number\subsubsecno) }#1}

\def\@th#1{\global\def\currenvir{th}\medbreak\noindent
    \global\advance\subsecno by 1 
    {\bf #1 \the\secno.\the\subsecno.--- }\begingroup\it}

\def\@@th{\endgroup\medbreak}

\def\@thwem#1{\global\def\currenvir{th}\medbreak\noindent
    \global\advance\subsecno by 1 
    {\bf #1 \the\secno.\the\subsecno.--- }\begingroup}

\def\@@thwem{\endgroup\medbreak}

\def\@thwn#1{\medbreak\noindent{\bf #1 .---}\begingroup\it}
\def\@@thwn{\endgroup\medbreak}
\def\@formula{\global\def\currenvir{formula}\global\advance\subsubsecno by 1 $$}

\def\@@formula{\leqno(\the\secno.\the\subsecno.\the\subsubsecno)$$}

\long\def\comment#1{\relax}

\let\efl\hfl

\let\sfl\vfl
\def\nfl#1#2{\uTo^{#1}_{#2}}

\def\@diagram{\diagram[moreoptions]}
\def\@@diagram{\enddiagram}
\def\mono{\rInto^{}_{}}
\def\iso{\,\smash{\mathop{\rightarrow}\limits^{\sim}}\,}

\def\lra{\rTo^{}_{}}
\diagramstyle{small,midshaft}

\catcode`\@=12

\catcode`\@=11

\def\undefine#1{\let#1\undefined}
\def\newsymbol#1#2#3#4#5{\let\next@\relax
 \ifnum#2=\@ne\let\next@\msafam@\else
 \ifnum#2=\tw@\let\next@\msbfam@\fi\fi
 \mathchardef#1="#3\next@#4#5}
\def\mathhexbox@#1#2#3{\relax
 \ifmmode\mathpalette{}{\m@th\mathchar"#1#2#3}%
 \else\leavevmode\hbox{$\m@th\mathchar"#1#2#3$}\fi}
\def\hexnumber@#1{\ifcase#1 0\or 1\or 2\or 3\or 4\or 5\or 6\or 7\or 8\or
 9\or A\or B\or C\or D\or E\or F\fi}

\font\tenmsa=msam10
\font\sevenmsa=msam7
\font\fivemsa=msam5
\newfam\msafam
\textfont\msafam=\tenmsa
\scriptfont\msafam=\sevenmsa
\scriptscriptfont\msafam=\fivemsa
\edef\msafam@{\hexnumber@\msafam}
\mathchardef\dabar@"0\msafam@39
\def\dashrightarrow{\mathrel{\dabar@\dabar@\mathchar"0\msafam@4B}}
\def\dashleftarrow{\mathrel{\mathchar"0\msafam@4C\dabar@\dabar@}}

\def\ulcorner{\delimiter"4\msafam@70\msafam@70 }
\def\urcorner{\delimiter"5\msafam@71\msafam@71 }
\def\llcorner{\delimiter"4\msafam@78\msafam@78 }
\def\lrcorner{\delimiter"5\msafam@79\msafam@79 }
\def\yen{{\mathhexbox@\msafam@55 }}
\def\checkmark{{\mathhexbox@\msafam@58 }}
\def\circledR{{\mathhexbox@\msafam@72 }}
\def\maltese{{\mathhexbox@\msafam@7A }}

\font\tenmsb=msbm10
\font\sevenmsb=msbm7
\font\fivemsb=msbm5
\newfam\msbfam
\textfont\msbfam=\tenmsb
\scriptfont\msbfam=\sevenmsb
\scriptscriptfont\msbfam=\fivemsb
\edef\msbfam@{\hexnumber@\msbfam}
\def\Bbb#1{{\fam\msbfam\relax#1}}
\def\widehat#1{\setbox\z@\hbox{$\m@th#1$}%
 \ifdim\wd\z@>\tw@ em\mathaccent"0\msbfam@5B{#1}%
 \else\mathaccent"0362{#1}\fi}
\def\widetilde#1{\setbox\z@\hbox{$\m@th#1$}%
 \ifdim\wd\z@>\tw@ em\mathaccent"0\msbfam@5D{#1}%
 \else\mathaccent"0365{#1}\fi}

\font\teneufm=eufm10
\font\seveneufm=eufm7
\font\fiveeufm=eufm5
\newfam\eufmfam
\textfont\eufmfam=\teneufm
\scriptfont\eufmfam=\seveneufm
\scriptscriptfont\eufmfam=\fiveeufm

\let\goth\frak

\newsymbol\boxdot 1200
\newsymbol\boxplus 1201
\newsymbol\boxtimes 1202
\newsymbol\square 1003
\newsymbol\blacksquare 1004
\newsymbol\centerdot 1205
\newsymbol\lozenge 1006
\newsymbol\blacklozenge 1007
\newsymbol\circlearrowright 1308
\newsymbol\circlearrowleft 1309
\undefine\rightleftharpoons
\newsymbol\rightleftharpoons 130A
\newsymbol\leftrightharpoons 130B
\newsymbol\boxminus 120C
\newsymbol\Vdash 130D
\newsymbol\Vvdash 130E
\newsymbol\vDash 130F
\newsymbol\twoheadrightarrow 1310
\newsymbol\twoheadleftarrow 1311
\newsymbol\leftleftarrows 1312
\newsymbol\rightrightarrows 1313
\newsymbol\upuparrows 1314
\newsymbol\downdownarrows 1315
\newsymbol\upharpoonright 1316
 \let\restriction\upharpoonright
\newsymbol\downharpoonright 1317
\newsymbol\upharpoonleft 1318
\newsymbol\downharpoonleft 1319
\newsymbol\rightarrowtail 131A
\newsymbol\leftarrowtail 131B
\newsymbol\leftrightarrows 131C
\newsymbol\rightleftarrows 131D
\newsymbol\Lsh 131E
\newsymbol\Rsh 131F
\newsymbol\rightsquigarrow 1320
\newsymbol\leftrightsquigarrow 1321
\newsymbol\looparrowleft 1322
\newsymbol\looparrowright 1323
\newsymbol\circeq 1324
\newsymbol\succsim 1325
\newsymbol\gtrsim 1326
\newsymbol\gtrapprox 1327
\newsymbol\multimap 1328
\newsymbol\therefore 1329
\newsymbol\because 132A
\newsymbol\doteqdot 132B
 
\newsymbol\triangleq 132C
\newsymbol\precsim 132D
\newsymbol\lesssim 132E
\newsymbol\lessapprox 132F
\newsymbol\eqslantless 1330
\newsymbol\eqslantgtr 1331
\newsymbol\curlyeqprec 1332
\newsymbol\curlyeqsucc 1333
\newsymbol\preccurlyeq 1334
\newsymbol\leqq 1335
\newsymbol\leqslant 1336
\newsymbol\lessgtr 1337
\newsymbol\backprime 1038
\newsymbol\risingdotseq 133A
\newsymbol\fallingdotseq 133B
\newsymbol\succcurlyeq 133C
\newsymbol\geqq 133D
\newsymbol\geqslant 133E
\newsymbol\gtrless 133F
\newsymbol\sqsubset 1340
\newsymbol\sqsupset 1341
\newsymbol\vartriangleright 1342
\newsymbol\vartriangleleft 1343
\newsymbol\trianglerighteq 1344
\newsymbol\trianglelefteq 1345
\newsymbol\bigstar 1046
\newsymbol\between 1347
\newsymbol\blacktriangledown 1048
\newsymbol\blacktriangleright 1349
\newsymbol\blacktriangleleft 134A
\newsymbol\vartriangle 134D
\newsymbol\blacktriangle 104E
\newsymbol\triangledown 104F
\newsymbol\eqcirc 1350
\newsymbol\lesseqgtr 1351
\newsymbol\gtreqless 1352
\newsymbol\lesseqqgtr 1353
\newsymbol\gtreqqless 1354
\newsymbol\Rrightarrow 1356
\newsymbol\Lleftarrow 1357
\newsymbol\veebar 1259
\newsymbol\barwedge 125A
\newsymbol\doublebarwedge 125B
\undefine\angle
\newsymbol\angle 105C
\newsymbol\measuredangle 105D
\newsymbol\sphericalangle 105E
\newsymbol\varpropto 135F
\newsymbol\smallsmile 1360
\newsymbol\smallfrown 1361
\newsymbol\Subset 1362
\newsymbol\Supset 1363
\newsymbol\Cup 1264
 
\newsymbol\Cap 1265
 
\newsymbol\curlywedge 1266
\newsymbol\curlyvee 1267
\newsymbol\leftthreetimes 1268
\newsymbol\rightthreetimes 1269
\newsymbol\subseteqq 136A
\newsymbol\supseteqq 136B
\newsymbol\bumpeq 136C
\newsymbol\Bumpeq 136D
\newsymbol\lll 136E
 
\newsymbol\ggg 136F
 
\newsymbol\circledS 1073
\newsymbol\pitchfork 1374
\newsymbol\dotplus 1275
\newsymbol\backsim 1376
\newsymbol\backsimeq 1377
\newsymbol\complement 107B
\newsymbol\intercal 127C
\newsymbol\circledcirc 127D
\newsymbol\circledast 127E
\newsymbol\circleddash 127F
\newsymbol\lvertneqq 2300
\newsymbol\gvertneqq 2301
\newsymbol\nleq 2302
\newsymbol\ngeq 2303
\newsymbol\nless 2304
\newsymbol\ngtr 2305
\newsymbol\nprec 2306
\newsymbol\nsucc 2307
\newsymbol\lneqq 2308
\newsymbol\gneqq 2309
\newsymbol\nleqslant 230A
\newsymbol\ngeqslant 230B
\newsymbol\lneq 230C
\newsymbol\gneq 230D
\newsymbol\npreceq 230E
\newsymbol\nsucceq 230F
\newsymbol\precnsim 2310
\newsymbol\succnsim 2311
\newsymbol\lnsim 2312
\newsymbol\gnsim 2313
\newsymbol\nleqq 2314
\newsymbol\ngeqq 2315
\newsymbol\precneqq 2316
\newsymbol\succneqq 2317
\newsymbol\precnapprox 2318
\newsymbol\succnapprox 2319
\newsymbol\lnapprox 231A
\newsymbol\gnapprox 231B
\newsymbol\nsim 231C
\newsymbol\ncong 231D
\newsymbol\diagup 231E
\newsymbol\diagdown 231F
\newsymbol\varsubsetneq 2320
\newsymbol\varsupsetneq 2321
\newsymbol\nsubseteqq 2322
\newsymbol\nsupseteqq 2323
\newsymbol\subsetneqq 2324
\newsymbol\supsetneqq 2325
\newsymbol\varsubsetneqq 2326
\newsymbol\varsupsetneqq 2327
\newsymbol\subsetneq 2328
\newsymbol\supsetneq 2329
\newsymbol\nsubseteq 232A
\newsymbol\nsupseteq 232B
\newsymbol\nparallel 232C
\newsymbol\nmid 232D
\newsymbol\nshortmid 232E
\newsymbol\nshortparallel 232F
\newsymbol\nvdash 2330
\newsymbol\nVdash 2331
\newsymbol\nvDash 2332
\newsymbol\nVDash 2333
\newsymbol\ntrianglerighteq 2334
\newsymbol\ntrianglelefteq 2335
\newsymbol\ntriangleleft 2336
\newsymbol\ntriangleright 2337
\newsymbol\nleftarrow 2338
\newsymbol\nrightarrow 2339
\newsymbol\nLeftarrow 233A
\newsymbol\nRightarrow 233B
\newsymbol\nLeftrightarrow 233C
\newsymbol\nleftrightarrow 233D
\newsymbol\divideontimes 223E
\newsymbol\varnothing 203F
\newsymbol\nexists 2040
\newsymbol\Finv 2060
\newsymbol\Game 2061
\newsymbol\mho 2066
\newsymbol\eth 2067
\newsymbol\eqsim 2368
\newsymbol\beth 2069
\newsymbol\gimel 206A
\newsymbol\daleth 206B
\newsymbol\lessdot 236C
\newsymbol\gtrdot 236D
\newsymbol\ltimes 226E
\newsymbol\rtimes 226F
\newsymbol\shortmid 2370
\newsymbol\shortparallel 2371
\newsymbol\smallsetminus 2272
\newsymbol\thicksim 2373
\newsymbol\thickapprox 2374
\newsymbol\approxeq 2375
\newsymbol\succapprox 2376
\newsymbol\precapprox 2377
\newsymbol\curvearrowleft 2378
\newsymbol\curvearrowright 2379
\newsymbol\digamma 207A
\newsymbol\varkappa 207B
\newsymbol\Bbbk 207C
\newsymbol\hslash 207D
\undefine\hbar
\newsymbol\hbar 207E
\newsymbol\backepsilon 237F

\catcode`\@=12

\font\tensc=cmcsc10 % text small capitals
\newfam\scfam \def\sc{\fam\scfam\tensc} % \sc is family 8
\textfont\scfam=\tensc

\font\sevenit=cmti7
\let\sevensl=\sevenit        \let\seventt=\sevenrm
 
\let\fourrm=\fiverm          \let\fouri=\fivei   \let\foursy=\fivesy
\let\fourbf=\fivebf
 
\ifx\sevenpoint\undefined
   \def\sevenpoint{\def\rm{\fam0\sevenrm} %switch to a 7-point type
       \textfont0=\sevenrm  \scriptfont0=\fiverm \scriptscriptfont0=\fourrm
       \textfont1=\seveni   \scriptfont1=\fivei  \scriptscriptfont1=\fouri
       \textfont2=\sevensy  \scriptfont2=\fivesy \scriptscriptfont2=\foursy
       \textfont3=\tenex    \scriptfont3=\tenex  \scriptscriptfont3=\tenex
       \textfont\itfam=\sevenit   \def\it{\fam\itfam\sevenit} \let\em\it%
       \textfont\slfam=\sevensl   \def\sl{\fam\slfam\sevensl}%
       \textfont\ttfam=\seventt   \def\tt{\fam\ttfam\seventt}%
       \textfont\bffam=\sevenbf   \scriptfont\bffam=\fivebf
        \scriptscriptfont\bffam=\fourbf   \def\bf{\fam\bffam\sevenbf}%
       \normalbaselineskip=10pt
       \setbox\strutbox=\hbox{\vrule height7.5pt depth2.5pt width0pt}%
       \normalbaselines\rm}
   \fi
\font\eightrm=cmr8           \font\sixrm=cmr6
\font\eighti=cmmi8           \font\sixi=cmmi6
\font\eightsy=cmsy8          \font\sixsy=cmsy6
\font\eightbf=cmbx8          \font\sixbf=cmbx6
\font\eightsl=cmsl8          \font\eighttt=cmtt8    \font\eightit=cmti8
\font\eightcsc=cmcsc10 scaled 800
\ifx\eightpoint\undefined
   \def\eightpoint{\def\rm{\fam0\eightrm} %switch to a 8-point type
       \textfont0=\eightrm  \scriptfont0=\sixrm  \scriptscriptfont0=\fiverm
       \textfont1=\eighti   \scriptfont1=\sixi   \scriptscriptfont1=\fivei
       \textfont2=\eightsy  \scriptfont2=\sixsy  \scriptscriptfont2=\fivesy
       \textfont3=\tenex    \scriptfont3=\tenex  \scriptscriptfont3=\tenex
       \textfont\itfam=\eightit   \def\it{\fam\itfam\eightit} \let\em\it%
       \textfont\slfam=\eightsl   \def\sl{\fam\slfam\eightsl}%
       \textfont\ttfam=\eighttt   \def\tt{\fam\ttfam\eighttt}%
       \textfont\bffam=\eightbf   \scriptfont\bffam=\sixbf
        \scriptscriptfont\bffam=\fivebf   \def\bf{\fam\bffam\eightbf}%
       \textfont\scfam=\eightcsc  \def\sc{\fam\scfam\eightcsc}%
       \normalbaselineskip=10pt
       \setbox\strutbox=\hbox{\vrule height7.5pt depth2.5pt width0pt}%
       \normalbaselines\rm}
   \fi

\font\ninerm=cmr9            \font\sixrm=cmr6
\font\ninei=cmmi9            \font\sixi=cmmi6
\font\ninesy=cmsy9           \font\sixsy=cmsy6
\font\ninebf=cmbx9           \font\sixbf=cmbx6
\font\ninesl=cmsl9           \font\ninett=cmtt9      \font\nineit=cmti9
\font\ninecsc=cmcsc10 scaled 900
 
\ifx\ninepoint\undefined
   \def\ninepoint{\def\rm{\fam0\ninerm} %switch to a 9-point type
       \textfont0=\ninerm  \scriptfont0=\sixrm  \scriptscriptfont0=\fiverm
       \textfont1=\ninei   \scriptfont1=\sixi   \scriptscriptfont1=\fivei
       \textfont2=\ninesy  \scriptfont2=\sixsy  \scriptscriptfont2=\fivesy
       \textfont3=\tenex   \scriptfont3=\tenex  \scriptscriptfont3=\tenex
       \textfont\itfam=\nineit   \def\it{\fam\itfam\nineit} \let\em\it%
       \textfont\slfam=\ninesl   \def\sl{\fam\slfam\ninesl}%
       \textfont\ttfam=\ninett   \def\tt{\fam\ttfam\ninett}%
       \textfont\bffam=\ninebf   \scriptfont\bffam=\sixbf
        \scriptscriptfont\bffam=\fivebf   \def\bf{\fam\bffam\ninebf}%
       \textfont\scfam=\ninecsc  \def\sc{\fam\scfam\ninecsc}%
       \normalbaselineskip=11pt
       \setbox\strutbox=\hbox{\vrule height8pt depth3pt width0pt}%
       \normalbaselines\rm}
   \fi

\def\tenpoint{\def\rm{\fam0\tenrm}% switch to 10-point type
    \textfont0=\tenrm  \scriptfont0=\sevenrm  \scriptscriptfont0=\fiverm
    \textfont1=\teni   \scriptfont1=\seveni   \scriptscriptfont1=\fivei
    \textfont2=\tensy  \scriptfont2=\sevensy  \scriptscriptfont2=\fivesy
    \textfont3=\tenex  \scriptfont3=\tenex    \scriptscriptfont3=\tenex
    \textfont\itfam=\tenit   \def\it{\fam\itfam\tenit} \let\em\it%
    \textfont\slfam=\tensl   \def\sl{\fam\slfam\tensl}%
    \textfont\ttfam=\tentt   \def\tt{\fam\ttfam\tentt}%
    \textfont\bffam=\tenbf   \scriptfont\bffam=\sevenbf
    \scriptscriptfont\bffam=\fivebf  \def\bf{\fam\bffam\tenbf}%
    \textfont\scfam=\tencsc  \def\sc{\fam\scfam\tencsc}%
    \normalbaselineskip=12pt
    \setbox\strutbox=\hbox{\vrule height8.5pt depth 3.5pt width0pt}%
    \normalbaselines\rm}

\font\elevenrm=cmr10    scaled \magstephalf
\font\eleveni=cmmi10    scaled \magstephalf
\font\elevensy=cmsy10   scaled \magstephalf
\font\elevenex=cmex10   scaled \magstephalf
\font\elevenbf=cmbx10   scaled \magstephalf
\font\elevensl=cmsl10   scaled \magstephalf
\font\eleventt=cmtt10   scaled \magstephalf
\font\elevenit=cmti10   scaled \magstephalf
\font\elevencsc=cmcsc10 scaled \magstephalf
 
\font\eightrm=cmr8           \font\sixrm=cmr6
\font\eighti=cmmi8           \font\sixi=cmmi6
\font\eightsy=cmsy8          \font\sixsy=cmsy6
\font\eightbf=cmbx8          \font\sixbf=cmbx6
 
\ifx\elevenpoint\undefined
   \def\elevenpoint{\def\rm{\fam0\elevenrm}% switch to 11-point type
       \textfont0=\elevenrm \scriptfont0=\eightrm \scriptscriptfont0=\sixrm
       \textfont1=\eleveni  \scriptfont1=\eighti  \scriptscriptfont1=\sixi
       \textfont2=\elevensy \scriptfont2=\eightsy \scriptscriptfont2=\sixsy
       \textfont3=\elevenex \scriptfont3=\elevenex\scriptscriptfont3=\elevenex
       \textfont\itfam=\elevenit  \def\it{\fam\itfam\elevenit} \let\em\it%
       \textfont\slfam=\elevensl  \def\sl{\fam\slfam\elevensl}%
       \textfont\ttfam=\eleventt  \def\tt{\fam\ttfam\eleventt}%
       \textfont\bffam=\elevenbf  \scriptfont\bffam=\eightbf
        \scriptscriptfont\bffam=\sixbf  \def\bf{\fam\bffam\elevenbf}%
       \textfont\scfam=\elevencsc \def\sc{\fam\scfam\elevencsc}%
       \normalbaselineskip=14pt
       \setbox\strutbox=\hbox{\vrule height9pt depth4pt width0pt}%
       \normalbaselines\rm}
   \fi
\font\twelverm=cmr12
\font\twelvei=cmmi12
\font\twelvesy=cmsy10 scaled \magstep1
\font\twelveex=cmex10 scaled \magstep1
\font\twelvebf=cmbx12
\font\twelvesl=cmsl12
\font\twelvett=cmtt12
\font\twelveit=cmti12
\font\twelvecsc=cmcsc10 scaled \magstep1
 
\font\ninerm=cmr9            \font\sevenrm=cmr7
\font\ninei=cmmi9            \font\seveni=cmmi7
\font\ninesy=cmsy9           \font\sevensy=cmsy7
\font\ninebf=cmbx9           \font\sevenbf=cmbx7
 
\ifx\twelvepoint\undefined
   \def\twelvepoint{\def\rm{\fam0\twelverm}% switch to 11-point type
       \textfont0=\twelverm \scriptfont0=\ninerm \scriptscriptfont0=\sevenrm
       \textfont1=\twelvei  \scriptfont1=\ninei  \scriptscriptfont1=\seveni
       \textfont2=\twelvesy \scriptfont2=\ninesy \scriptscriptfont2=\sevensy
       \textfont3=\twelveex \scriptfont3=\twelveex\scriptscriptfont3=\twelveex
       \textfont\itfam=\twelveit  \def\it{\fam\itfam\twelveit} \let\em\it%
       \textfont\slfam=\twelvesl  \def\sl{\fam\slfam\twelvesl}%
       \textfont\ttfam=\twelvett  \def\tt{\fam\ttfam\twelvett}%
       \textfont\bffam=\twelvebf  \scriptfont\bffam=\ninebf %
       \scriptscriptfont\bffam=\sevenbf  \def\bf{\fam\bffam\twelvebf}%
       \textfont\scfam=\twelvecsc \def\sc{\fam\scfam\twelvecsc}%
       \normalbaselineskip=14pt
       \setbox\strutbox=\hbox{\vrule height9.5pt depth4.5pt width0pt}%
       \normalbaselines\rm}
   \fi

\font\fourteenrm=cmr12 scaled \magstep1
\font\fourteeni=cmmi12 scaled \magstep1
\font\fourteensy=cmsy10 scaled \magstep2
\font\fourteenex=cmex10 scaled \magstep2
\font\fourteenbf=cmbx12 scaled \magstep1
\font\fourteensl=cmsl12 scaled \magstep1
\font\fourteentt=cmtt12 scaled \magstep1
\font\fourteenit=cmti12 scaled \magstep1
\font\fourteencsc=cmcsc10 scaled \magstep2
 
\ifx\fourteenpoint\undefined
   \def\fourteenpoint{\def\rm{\fam0\fourteenrm}% switch to 14-point type
       \textfont0=\fourteenrm \scriptfont0=\tenrm \scriptscriptfont0=\sevenrm
       \textfont1=\fourteeni  \scriptfont1=\teni  \scriptscriptfont1=\seveni
       \textfont2=\fourteensy \scriptfont2=\tensy \scriptscriptfont2=\sevensy
       \textfont3=\fourteenex \scriptfont3=\fourteenex
                              \scriptscriptfont3=\fourteenex
       \textfont\itfam=\fourteenit  \def\it{\fam\itfam\fourteenit} \let\em\it%
       \textfont\slfam=\fourteensl  \def\sl{\fam\slfam\fourteensl}%
       \textfont\ttfam=\fourteentt  \def\tt{\fam\ttfam\fourteentt}%
       \textfont\bffam=\fourteenbf  \scriptfont\bffam=\tenbf
        \scriptscriptfont\bffam=\sevenbf  \def\bf{\fam\bffam\fourteenbf}%
       \textfont\scfam=\fourteencsc \def\sc{\fam\scfam\fourteencsc}%
       \normalbaselineskip=17pt
       \setbox\strutbox=\hbox{\vrule height11.9pt depth6.3pt width0pt}%
       \normalbaselines\rm}
   \fi

\def\comp{{\Bbb C}}          %les nombres complexes
\def\proj{{\Bbb P}}          %l'espace projectif
\def\reln{{\Bbb Z}}          %les entiers relatifs
 
\def\droitep{\proj^{1}}     %droiteprojectif

\def\cf{{\it cf.\/}\ }
\def\ie{{\it i.e.\/}\ }

\def\nl{\hfill\break}\let\\\nl
\def\np{\vfill\eject}
\def\text#1{\leavevmode\hbox{#1}}

\newbox\boxa
\newbox\boxb
\newdimen\dima
\newdimen\dimb

\def\boxit#1#2{\setbox\boxa=\hbox{\kern#1{#2}\kern#1}%
  \dimen\dima=\ht\boxa \advance\dimen\dima by #1
  \dimen\dimb=\dp\boxa \advance\dimen\dimb by #1
  \setbox\boxa=\hbox{\vrule\box\boxa\vrule}
  \setbox\boxa=\vbox{\hrule\box\boxa\hrule}
  \advance\dimen\dima by 0.4pt \ht\boxa=\dimen\dima
  \advance\dimen\dimb by 0.4pt \dp\boxa=\dimen\dimb
  \box\boxa\relax}

\mathcode`A="7041 \mathcode`B="7042 \mathcode`C="7043
\mathcode`D="7044 \mathcode`E="7045 \mathcode`F="7046
\mathcode`G="7047 \mathcode`H="7048 \mathcode`I="7049
\mathcode`J="704A \mathcode`K="704B \mathcode`L="704C
\mathcode`M="704D \mathcode`N="704E \mathcode`O="704F
\mathcode`P="7050 \mathcode`Q="7051 \mathcode`R="7052
\mathcode`S="7053 \mathcode`T="7054 \mathcode`U="7055
\mathcode`V="7056 \mathcode`W="7057 \mathcode`X="7058
\mathcode`Y="7059 \mathcode`Z="705A

\def\frac#1#2{{#1\over #2}}
\def\cqfd{\kern 2truemm\unskip\penalty 500\vrule height 4pt depth 0pt
width 4pt\medbreak}

\def\no{n\up{o}\kern 2pt}
\def\moins{\mathrel{\hbox{\vrule height 3pt depth -2pt width 6pt}}}
\def\rond{\kern 1pt{\scriptstyle\circ}\kern 1pt}
\let\bk\backslash
\let\ra\rightarrow
\let\ol\overline
\let\ul\underline
\let\ort\perp
\def\osum{\mathop{\oplus}\limits}
\def\tvi{\vrule height 12pt depth 5pt width 0pt}
\def\tv{\tvi\vrule}
\def\cc#1{\hfill\kern .7em#1\kern .7em\hfill}

\def\End{\mathop{\rm End}\nolimits}
\def\Hom{\mathop{\rm Hom}\nolimits}

\def\Sym{\mathop{\rm Sym}\nolimits}
\def\ext{\mathop{\rm ext}\nolimits}

\def\Ext{\mathop{\rm Ext}\nolimits}
\def\Aut{\mathop{\rm Aut}\nolimits}
\def\Im{\mathop{\rm Im}\nolimits}
\def\Id{\mathop{\rm Id}\nolimits}
\def\Ker{\mathop{\rm Ker}\nolimits}
\def\Coker{\mathop{\rm Coker}}
\def\Spec{\mathop{\rm Spec}}

\def\det{\mathop{\rm det}\nolimits}
\def\Pic{\mathop{\rm Pic}\nolimits}
\def\Cl{\mathop{\rm Cl}\nolimits}

\def\dim{\mathop{\rm dim}\nolimits}
\def\codim{\mathop{\rm codim}\nolimits}

\def\Tr{\mathop{\rm Tr}\nolimits}

\def\Ad{\mathop{\rm Ad}\nolimits}
\def\Res{\mathop{\rm Res}\nolimits}
\def\Lie{\mathop{\rm Lie}\nolimits}
\def\Spin{\mathop{\rm Spin}\nolimits}

\def\limind{\mathop{\oalign{lim\cr\hidewidth$\longrightarrow
$\hidewidth\cr}}}
\def\limproj{\mathop{\oalign{lim\cr\hidewidth$\longleftarrow
$\hidewidth\cr}}}
\def\epi{\rightarrow \kern -3mm\rightarrow }

\newarrow{SL}-----
\newarrow{DL}=====
\newarrow{TL}33333

\def\myddar#1.#2.#3.#4.#5.#6{
\begin{diagram}[shortfall=-0.8pt,width=#6em,height=0.8em]
\circ&\rSL&\circ&\rSL&\circ&\rSL&\ \cdot\cdot\cdot\ &
\rSL&\circ&\rSL&\circ\cr
\hbox to 0pt{\hss$\scriptstyle#1$\hss}&&\hbox to 0pt{\hss$\scriptstyle#2$\hss}&&
\hbox to 0pt{\hss$\scriptstyle#3$\hss}&&&&
\hbox to 0pt{\hss$\scriptstyle#4$\hss}&&\hbox to 0pt{\hss$\scriptstyle#5$\hss}\cr
\end{diagram}}

\def\ddar#1.#2.#3.#4.#5.{\myddar#1.#2.#3.#4.#5.1}
\def\bigddar#1.#2.#3.#4.#5.{\myddar#1.#2.#3.#4.#5.2}

\def\myddbr#1.#2.#3.#4.#5.#6{
\tmpdimension=3.5pt\multiply\tmpdimension by #6
\ifdim\tmpdimension=7pt\tmpdimension=14pt\fi
\begin{diagram}[shortfall=-0.8pt,width=#6em,height=0.8em]
\circ&\rSL&\circ&\rSL&\circ&\rSL&\ \cdot\cdot\cdot\ &
\rSL&\circ&\rlap{\hskip\tmpdimension$>$}\rDL&\bullet\cr
\hbox to 0pt{\hss$\scriptstyle#1$\hss}&&\hbox to 0pt{\hss$\scriptstyle#2$\hss}&&
\hbox to 0pt{\hss$\scriptstyle#3$\hss}&&&&
\hbox to 0pt{\hss$\scriptstyle#4$\hss}&&\hbox to 0pt{\hss$\scriptstyle#5$\hss}\cr
\end{diagram}}

\def\ddbr#1.#2.#3.#4.#5.{\myddbr#1.#2.#3.#4.#5.1}
\def\bigddbr#1.#2.#3.#4.#5.{\myddbr#1.#2.#3.#4.#5.2}

\def\myddcr#1.#2.#3.#4.#5.#6{
\tmpdimension=3.5pt\multiply\tmpdimension by #6
\ifdim\tmpdimension=7pt\tmpdimension=14pt\fi
\begin{diagram}[shortfall=-0.8pt,width=#6em,height=0.8em]
\bullet&\rSL&\bullet&\rSL&\bullet&\rSL&\ \cdot\cdot\cdot\ &
\rSL&\bullet&\rlap{\hskip\tmpdimension$<$}\rDL&\circ\cr
\hbox to 0pt{\hss$\scriptstyle#1$\hss}&&\hbox to 0pt{\hss$\scriptstyle#2$\hss}&&
\hbox to 0pt{\hss$\scriptstyle#3$\hss}&&&&
\hbox to 0pt{\hss$\scriptstyle#4$\hss}&&\hbox to 0pt{\hss$\scriptstyle#5$\hss}\cr
\end{diagram}}

\def\ddcr#1.#2.#3.#4.#5.{\myddcr#1.#2.#3.#4.#5.1}
\def\bigddcr#1.#2.#3.#4.#5.{\myddcr#1.#2.#3.#4.#5.2}

\def\dddr#1.#2.#3.#4.#5.#6.{
\begin{diagram}[shortfall=-0.8pt,width=1em,height=0.8em,silent]
&&&&&&&&&\hskip8pt\hbox to 0pt{\hss$\scriptstyle#6$\hss}\cr
&&&&&&&&&\hskip8pt\circ\cr
\circ&\rSL&\circ&\rSL&\circ&\rSL&\ \cdot\cdot\cdot\ &
\rSL&\circ\kern-.5pt\ruSL\cr
\hbox to 0pt{\hss$\scriptstyle#1$\hss}&&\hbox to 0pt{\hss$\scriptstyle#2$\hss}&&
\hbox to 0pt{\hss$\scriptstyle#3$\hss}&&&&
\hbox to 0pt{\hss$\scriptstyle#4$\hss}\kern4pt\rdSL&\hskip8pt\circ\cr
&&&&&&&&&\hskip8pt\hbox to 0pt{\hss$\scriptstyle#5$\hss}\cr
\end{diagram}}

\def\bigdddr#1.#2.#3.#4.#5.#6.{
\begin{diagram}[shortfall=-0.8pt,width=2em,height=0.8em,silent]
&&&&&&&&&\hskip8pt\hbox to 0pt{\hss$\scriptstyle#6$\hss}\cr
&&&&&&&&&\hskip8pt\circ\cr
&&&&&&&&\kern24pt\ruSL[h=1.5em]\cr
\circ&\rSL&\circ&\rSL&\circ&\rSL&\ \cdot\cdot\cdot\ &
\rSL&\circ\cr
\hbox to 0pt{\hss$\scriptstyle#1$\hss}&&\hbox to 0pt{\hss$\scriptstyle#2$\hss}&&
\hbox to 0pt{\hss$\scriptstyle#3$\hss}&&&&
\hbox to 0pt{\hss$\scriptstyle#4$\hss}&\kern-16pt\rdSL[h=1.5em]\cr
&&&&&&&&&\hskip8pt\circ\cr
&&&&&&&&&\hskip8pt\hbox to 0pt{\hss$\scriptstyle#5$\hss}\cr
\end{diagram}}

\def\myddei#1.#2.#3.#4.#5.#6.#7{
\begin{diagram}[shortfall=-0.8pt,width=#7em,height=0.8em]
&&&&\circ&\hbox to 0pt{\hss$\scriptstyle#2$\hss}\cr
&&&&\dSL\cr
\circ&\rSL&\circ&\rSL&
\vbox{\hrule height 0.5pt width 0pt\hbox{$\circ$}}&\rSL&\circ&
\rSL&\circ\cr
\hbox to 0pt{\hss$\scriptstyle#1$\hss}&&\hbox to 0pt{\hss$\scriptstyle#3$\hss}&&
\hbox to 0pt{\hss$\scriptstyle#4$\hss}&&
\hbox to 0pt{\hss$\scriptstyle#5$\hss}&&\hbox to 0pt{\hss$\scriptstyle#6$\hss}\cr
\end{diagram}}

\def\ddei#1.#2.#3.#4.#5.#6.{\myddei#1.#2.#3.#4.#5.#6.1}
\def\bigddei#1.#2.#3.#4.#5.#6.{\myddei#1.#2.#3.#4.#5.#6.2}

\def\myddeii#1.#2.#3.#4.#5.#6.#7.#8{
\begin{diagram}[shortfall=-0.8pt,width=#8em,height=0.8em]
&&&&\circ&\hbox to 0pt{\hss$\scriptstyle#2$\hss}\cr
&&&&\uSL\cr
\circ&\rSL&\circ&\rSL&
\vbox{\hrule height 0.5pt width 0pt\hbox{$\circ$}}&\rSL&\circ&
\rSL&\circ&\rSL&\circ\cr
\hbox to 0pt{\hss$\scriptstyle#1$\hss}&&\hbox to 0pt{\hss$\scriptstyle#3$\hss}&&
\hbox to 0pt{\hss$\scriptstyle#4$\hss}&&\hbox to 0pt{\hss$\scriptstyle#5$\hss}&&
\hbox to 0pt{\hss$\scriptstyle#6$\hss}&&\hbox to 0pt{\hss$\scriptstyle#7$\hss}\cr
\end{diagram}}

\def\ddeii#1.#2.#3.#4.#5.#6.#7.{\myddeii#1.#2.#3.#4.#5.#6.#7.1}
\def\bigddeii#1.#2.#3.#4.#5.#6.#7.{\myddeii#1.#2.#3.#4.#5.#6.#7.2}

\def\myddeiii#1.#2.#3.#4.#5.#6.#7.#8.#9{
\begin{diagram}[shortfall=-0.8pt,width=#9em,height=0.8em]
&&&&\circ&\hbox to 0pt{\hss$\scriptstyle#2$\hss}\cr
&&&&\uSL\cr
\circ&\rSL&\circ&\rSL&
\vbox{\hrule height 0.5pt width 0pt\hbox{$\circ$}}&\rSL&\circ&
\rSL&\circ&\rSL&\circ&\rSL&\circ\cr
\hbox to 0pt{\hss$\scriptstyle#1$\hss}&&\hbox to 0pt{\hss$\scriptstyle#3$\hss}&&
\hbox to 0pt{\hss$\scriptstyle#4$\hss}&&\hbox to 0pt{\hss$\scriptstyle#5$\hss}&&
\hbox to 0pt{\hss$\scriptstyle#6$\hss}&&\hbox to 0pt{\hss$\scriptstyle#7$\hss}&&
\hbox to 0pt{\hss$\scriptstyle#8$\hss}\cr
\end{diagram}}

\def\ddeiii#1.#2.#3.#4.#5.#6.#7.#8.{\myddeiii#1.#2.#3.#4.#5.#6.#7.#8.1}
\def\bigddeiii#1.#2.#3.#4.#5.#6.#7.#8.{\myddeiii#1.#2.#3.#4.#5.#6.#7.#8.2}

\def\myddf#1.#2.#3.#4.#5{
\tmpdimension=3.5pt\multiply\tmpdimension by #5
\ifdim\tmpdimension=7pt\tmpdimension=14pt\fi
\begin{diagram}[shortfall=-0.8pt,width=#5em,height=0.8em]
\circ&\rSL&\circ&\rlap{\hskip\tmpdimension$>$}\rDL&\bullet&\rSL&\bullet\cr
\hbox to 0pt{\hss$\scriptstyle#1$\hss}&&
\hbox to 0pt{\hss$\scriptstyle#2$\hss}&&
\hbox to 0pt{\hss$\scriptstyle#3$\hss}&&
\hbox to 0pt{\hss$\scriptstyle#4$\hss}\cr
\end{diagram}}

\def\ddf#1.#2.#3.#4.{\myddf#1.#2.#3.#4.1}
\def\bigddf#1.#2.#3.#4.{\myddf#1.#2.#3.#4.2}

\def\ddg#1.#2.{
\tmpdimension=3.5pt
\begin{diagram}[shortfall=-0.8pt,width=1em,height=0.8em]
\rlap{\hskip 4.2pt$\vcenter{\hrule height 0.2pt depth 0.2pt width
16.5pt}$}
\circ
&\rlap{\hskip\tmpdimension$>$}\rDL
&\bullet\cr
\hbox to 0pt{\hss$\scriptstyle#1$\hss}&&
\hbox to 0pt{\hss$\scriptstyle#2$\hss}\cr
\end{diagram}}

\def\bigddg#1.#2.{
\tmpdimension=15pt
\begin{diagram}[shortfall=-0.8pt,width=2em,height=0.8em]
\rlap{\hskip 4.2pt$\vcenter{\hrule height 0.2pt depth 0.2pt width
36.5pt}$}
\circ
&\rlap{\hskip\tmpdimension$>$}\rDL
&\bullet\cr
\hbox to 0pt{\hss$\scriptstyle#1$\hss}&&
\hbox to 0pt{\hss$\scriptstyle#2$\hss}\cr
\end{diagram}}

{\catcode`\@=11
\gdef\AllowBreak{\global \dt@ptrue
\everycr {\noalign 
{\ifdt@p \global \dt@pfalse\else \penalty\interdisplaylinepenalty 
\fi}}}
}

\def\nc#1{\def#1}
\def\operatorname#1{\mathop{\rm #1}\nolimits}
\def\operatornamewithlimits#1{\mathop{\rm #1}}
\nc{\on}{\operatorname}
\nc{\onw}{\operatornamewithlimits}

\nc{\ra}{\rightarrow}
\nc{\GL}{\on{GL}}
\nc{\SL}{\on{SL}}
\nc{\U}{\on{U}}
\nc{\SO}{\on{SO}}
\nc{\Spin}{\on{Spin}}
\nc{\Sym}{\on{Sym}}
\nc{\Res}{\onw{Res}}
\nc{\Ext}{\on{Ext}}
\nc{\Tor}{\on{Tor}}
\nc{\Hom}{\on{Hom}}
\nc{\End}{\on{End}}
\nc{\Pic}{\on{Pic}}
\nc{\H}{\on{H}}
\nc{\Todd}{\on{Todd}}
\nc{\res}{\on{Res}}
\nc{\Spec}{\on{Spec}}
\nc{\Gr}{\on{Gr}}
\nc{\CoInd}{\on{CoInd}}
\nc{\Ind}{\on{Ind}}
\nc{\Vac}{\on{Vac}}
\nc{\Tr}{\on{Tr}}
\nc{\tensor}{\otimes}
\nc{\rang}{\on{rang}}
\nc{\ds}{\displaystyle}
\nc{\supp}{\on{supp}}
\nc{\Quot}{\on{Quot}}
\nc{\supps}{\supp_{s}}
\nc{\prof}{\on{prof}}
\nc{\dh}{\on{dh}}
\nc{\gr}{\on{gr}}
\nc{\id}{\on{id}}
\nc{\HyperExt}{\on{\tr{E}}xt}
\nc{\HyperH}{\on{\tr{H}}}
\nc{\Im}{\on{Im}}
\nc{\Ker}{\on{Ker}}
\nc{\Coker}{\on{Coker}}
\nc{\Aut}{\on{Aut}}
\nc{\codim}{codim}
\nc{\Cl}{\on{Cl}}
\nc{\Pf}{\on{Pf}}
\nc{\pf}{\on{pf}}
\nc{\ext}{\on{ext}}
\nc{\h}{\on{h}}
\nc{\mult}{\on{mult}}
\nc{\Groth}{\on{Groth}}
\nc{\Stab}{\on{Stab}}
\nc{\droitep}{\proj_{1}}
\nc{\ie}{{\it i.e. }}
\nc{\ul}{\underline}
\nc{\ol}{\overline}
\nc{\inject}{\mono}
\nc{\cf}{{\it cf.} }
\nc{\ort}{\perp}
\nc{\R}{\on{R}}
\nc{\tvi}{\vrule height 12pt depth 5pt width 0pt}
\nc{\tv}{\tvi\vrule}
\nc{\osum}{\onw{\oplus}}
\nc{\np}{\clearpage}
\nc{\Vir}{\on{Vir}}
\nc{\no}{\mathop{\raisebox{-2.5pt}{$\stackrel{\textstyle\circ}{\circ}$}}}
\nc{\ensfrac}#1#2{\raisebox{+4pt}{$#1$}/\raisebox{-4pt}{$#2$}}
\nc{\motimes}{\dot{\otimes}}
\nc{\moins}{\mathrel{\hbox{\vrule height 3pt depth -2pt width 6pt}}}
\nc{\moinss}{\mathrel{\hbox{\kern1pt\vrule height 2.3pt depth -1.6pt
width 4.2pt\kern1pt}}}
\nc{\ob}{\on{ob}}
\nc{\car}{\on{char}}
\nc{\Lie}{\on{Lie}}
\nc{\Isom}{\on{Isom}}
\nc{\Mor}{\on{Mor}}
\nc{\extern}{\onw{\boxtimes}}
\nc{\Ad}{\on{Ad}}
\nc{\rank}{\on{rank}}
\nc{\bk}{\backslash}
\nc{\lla}{\longleftarrow}
\nc{\limind}{\mathop{\oalign{lim\cr\hidewidth$\longrightarrow$\hidewidth\cr}}}
\nc{\limproj}{\mathop{\oalign{lim\cr\hidewidth$\longleftarrow$\hidewidth\cr}}}
\nc{\mono}{\hookrightarrow}
\nc{\omal}{\onw{\otimes}}
\nc{\Sp}{\on{Sp}}
\nc{\restriction}#1{_{|#1}}

\nc{\Pf}{\on{Pf}}
\nc{\det}{\on{det}}
\nc{\codim}{\on{codim}}

\nc{\Det}{\on{Det}}
\nc{\Id}{\on{Id}}
\nc{\ev}{\on{ev}}
\nc{\idbb}{{\Bbb I}}
\let\scr\cal

\nc{\g}{{\goth{g}}}
\nc{\gt}{\widetilde{\goth{g}}}
\nc{\Lg}{L\g}
\nc{\Lkg}{L^{k}\g}
\nc{\Lpolg}{L_{\rm{p}}\g}
\nc{\Lgmm}{L^{<0}\g}
\nc{\Lgpp}{L^{>0}\g}
\nc{\Lgh}{\widehat{\Lg}}
\nc{\Lkgh}{\widehat{\Lkg}}
\nc{\Lpolgh}{\widehat{\Lpolg}}
\nc{\Lgp}{L^{+}\g}
\nc{\Lgm}{L^{-}\g}
\nc{\ph}{\hat{\goth{p}}}
\nc{\Lgppol}{L^{+}_{_{pol}}\g}
\nc{\LPoneg}{L_{\droitep}\g}
\nc{\Lgph}{\widehat{\Lgp}}
\nc{\LgX}{L_{X}\g}
\nc{\LgU}{L_{_{U}}\g}
\nc{\LgUq}{L_{_{U-q}}\g}
\nc{\LGp}{L^{+}G}
\nc{\LGm}{L^{-}G}
\nc{\LGmm}{L^{<0}G}
\nc{\LG}{LG}
\nc{\LGh}{\widehat{\LG}}
\nc{\LGph}{\widehat{\LGp}}
\nc{\LGX}{L_{X}G}
\nc{\LGXi}{L_{X}^{i}G}
\nc{\LGXone}{L_{X}^{1}G}
\nc{\LGXip}{L_{X}^{i+1}G}
\nc{\LGXN}{\LGX(N)}
\nc{\LGXiN}{\LGXi(N)}
\nc{\LGN}{\LG(N)}
\nc{\LGzero}{\LG(0)}
\nc{\LSL}{LSL_{r}}
\nc{\LSLtwo}{LSL_{2}}
\nc{\LSLh}{\widehat{\LSL}}
\nc{\LSLp}{L^{+}SL_{r}}
\nc{\LSLX}{L_{X}SL_{r}}
\nc{\LSLN}{\LSL(N)}
\nc{\Lsl}{L{\goth{sl}}_{r}}
\nc{\LslV}{L{\goth{sl}}(V)}
\nc{\LslVh}{\widehat{\LslV}}
\nc{\Lslh}{\widehat{\Lsl}}
\nc{\rond}{\circ}
\nc{\G}{{\bf G}}
\nc{\HG}{{\cal{H}}}
\nc{\M}{{\cal{M}}_{G,X}}
\nc{\MEeight}{{\cal{M}}_{E_8,X}}
\nc{\MEseven}{{\cal{M}}_{E_7,X}}
\nc{\MEsix}{{\cal{M}}_{E_6,X}}
\nc{\MFfour}{{\cal{M}}_{F_4,X}}
\nc{\MSpin}{{\cal{M}}_{\Spin_{r},X}}
\nc{\MSO}{{\cal{M}}_{SO_{r},X}}
\nc{\MSLTF8}{{\cal{M}}_{SL_{248},X}}
\nc{\MGtwo}{{\cal{M}}_{G_2},X}
\nc{\MSOzero}{{\cal{M}}_{SO_{r},X}^{0}}
\nc{\MSOone}{{\cal{M}}_{SO_{r},X}^{1}}
\nc{\Mproj}{{\cal{M}}_{G,\proj^{1}}}
\nc{\MM}{{\cal{M}}_{G^{\prime}}}
\nc{\Q}{{\cal{Q}}_G}
\nc{\QEeight}{{\cal{Q}}_{E_8}}
\nc{\QEseven}{{\cal{Q}}_{E_7}}
\nc{\QEsix}{{\cal{Q}}_{E_6}}
\nc{\QFfour}{{\cal{Q}}_{F_4}}

\nc{\QSL}{{\cal{Q}}_{SL_{r}}}
\nc{\QSLtwo}{{\cal{Q}}_{SL_{2}}}
\nc{\QSLN}{{\cal{Q}}_{SL_{r}}(N)}
\nc{\QN}{{\cal{Q}}_G(N)}
\nc{\Qzero}{{\cal{Q}}_G(0)}
\nc{\QNred}{{\cal{Q}}_G(N)_{\on{red}}}
\nc{\Qp}{{\cal{Q}}_G^{p_{i+1}}}
\nc{\Qpar}{{\cal{Q}}^{par}_{G}(\ul{p},\ul{P})}
\nc{\Mpar}{{\cal{M}}^{par}_{G,X}(\ul{p},\ul{P})}
\nc{\ModSO}{M_{SO_{r},X}}
\nc{\ModSL}{M_{SL_{r},X}}
\nc{\ModSpin}{M_{\Spin_{r},X}}
\nc{\ModSpinEight}{M_{\Spin_{8},X}}
\nc{\ModEeight}{M_{E_{8},X}}
\nc{\ModGLzero}{M_{GL_{r},X}(0)}
\nc{\Modzero}{M_{G,X}(0)}
\nc{\ModSOreg}{M_{SO_{r},X}^{reg}}
\nc{\ModSOone}{M_{SO_{r},X}^{1}}
\nc{\ModSOzero}{M_{SO_{r},X}^{0}}
\nc{\Modt}{M_{G,X}^{\tau}}
\nc{\Mod}{M_{G,X}}
\nc{\Modreg}{M_{G,X}^{reg}}
\nc{\Modtreg}{M^{\tau,reg}_{G,X}}
\nc{\Mt}{{\cal{M}}_{G,X}^{\tau}}
\nc{\Rt}{R_{G,X}^{\tau}}
\nc{\Mtss}{{\cal{M}}^{\tau,ss}_{G,X}}
\nc{\Mss}{{\cal{M}}^{ss}_{G,X}}
\nc{\Mreg}{{\cal{M}}^{reg}_{G,X}}
\nc{\MSL}{{\cal{M}}_{SL_{r},X}}
\nc{\MGL}{{\cal{M}}_{GL_{r},X}}
\nc{\LGmmm}{L^{-}G_{m}}
\nc{\ma}{{\goth{M}}}
\nc{\Quad}{{\cal{R}}}
\nc{\Mum}{\on{Mum}}
\nc{\gpol}{\widehat\g}
\nc{\gpolp}{\widehat\g_{+}}
\nc{\gpolm}{\widehat\g_{-}}
\nc{\Cp}{{\goth{C}}^{+}}
\nc{\Cm}{{\goth{C}}^{-}}

\begin{document}
\baselineskip=15pt
%\null\bigskip
\centerline{\twelvepoint\bf On Moduli of G-bundles over Curves for
exceptional G}
\bigskip
\centerline{Christoph \sc{Sorger}}
\bigskip\bigskip

\section{Introduction}

\ind Let $G$ be a simple and simply connected complex Lie group, $\g$ its Lie
algebra. In the following, I remove the restriction ``$G$ is of classical type
or
$G_2$'' made on $G$ in the papers of Beauville, Laszlo and
myself \cite{L-S:verlinde},\cite{B-L-S:picard}
on the moduli of principal $G$-bundles on a curve.  As I will
just ``patch" the missing technical points, this note should
be seen as an appendix to the above cited papers.

\ind Let $\M$ be the stack of $G$-bundles on the smooth, connected
and projective algebraic curve $X$ of genus $g$. If
$\rho:G\ra\SL_{r}$ is a representation of $G$, denote by
${\cal{D}}_{\rho}$ the pullback of the determinant bundle
\cite{D-N:picard} under the morphism
$\M\ra\MSL$ defined by extension of the structure group.  Associate
to $G$ the number $d(G)$ and the representation $\rho(G)$ as follows:
$${\eightpoint\vbox{\offinterlineskip\def\hline{\noalign{\hrule}}
  \halign{\cc{$#$}&\tv#&\cc{$#$}&\tv#&\cc{$#$}&\tv#&%
          \cc{$#$}&\tv#&\cc{$#$}&\tv#&\cc{$#$}&\tv#&%
          \cc{$#$}&\tv#&\cc{$#$}&\tv#&\cc{$#$}&\tv#&%
          \cc{$#$}\cr
\text{Type of }G&&A_{r}&&B_{r}\,(r\geq 3)&&C_{r}&&D_{r}\, (r\geq
4)&&E_{6}&&E_{7}&&E_{8}&&F_{4}&&G_{2}\cr\hline
d(G)&&1&&2&&1&&2&&6&&12&&60&&6&&2\cr\hline
\rho(G)&&
\varpi_{1}&&\varpi_{1}&&\varpi_{1}&&\varpi_{1}&&\varpi_{6}
&&\varpi_{7}&&\varpi_{8}
&&\varpi_{4}&&\varpi_{1}\cr}}}$$

\begin{th}{Theorem}\label{th:pic} There is a line bundle ${\cal{L}}$
on $\M$ such that 
$\Pic(\M)\iso\reln{\cal{L}}$. Moreover we may choose ${\cal{L}}$ in
such a way that 
${\cal{L}}^{\otimes d(G)}={\cal{D}}_{\rho(G)}$.
\end{th} The above theorem is proved, for classical $G$ and $G_{2}$,
in \cite{L-S:verlinde} where it also shown that the space of
sections $H^{0}(\M,{\cal{L}}^{\ell})$ may be identified to the space
of conformal blocks $B_{G,X}(\ell;p;0)$ (see
(\ref{def:conf-blocks}) for its definition). Now, once the generator
of the Picard group is known in the exceptional cases, this
identification is also valid in general, as well what happens when
one considers additionally parabolic structures as we did in
\cite{L-S:verlinde} (theorems 1.1 and 1.2).

\ind In fact, as we will see, to prove theorem \ref{th:pic} for the
exceptional groups it is enough to prove the existence of the 60-th
root of
${\cal{D}}_{\varpi_{8}}$ on $\MEeight$. This will be deduced from
the splitting of a certain central extension, which in turn will
follow from the fact that $B_{E_8,X}(1;p;0)$ is one dimensional in
any genus $g$ as predicted by the Verlinde formula. However, in our particulier case
we don't need the  Verlinde formula in order to prove the last statement: it
will follow directly from the decomposition formulas. 
\ind Suppose $g(X)\geq 2$. For the coarse moduli spaces $\Mod$ of
semi-stable
$G$-bundles, we will see that the roots of the determinant bundle of
theorem
\ref{th:pic} do only exist on the open subset of regularly stable
$G$-bundles which, as shown in 
\cite{B-L-S:picard}, has as consequence the following:

\begin{th}{Theorem}\label{th:local_factoriality} Let $G$ be
semi-simple and $\tau\in\pi_{1}(G)$. Then $\Mod^{\tau}$ is locally factorial if and
only if
$G$ is special in the sens of Serre.
\end{th}

\ind Note that $\dim H^{0}(\MEeight,{\cal{L}})=\dim
B_{E_{8},X}(1;p;0)=1$  has the somehow surprising consequence that
the stack $\MEeight$ and (for
$g(X)\geq 2$) the normal variety $\ModEeight$ have a {\em canonical}
hypersurface.

\ind  
I would like to thank C. Teleman for pointing out that a reference I used
in a previous version of this paper was incomplete and mention his
preprint \cite{Te:BWB}, which contains a different,
topological approach to theorem \ref{th:pic}.

\section{Conformal Blocks}

\subsection{\em Affine Lie algebras. }\label{subsec:Lie-general-set-up} Let
$\g$ be a simple finite dimensional Lie algebra of rank $r$ over
$\comp$. Let $P$ be the weight lattice, $P_{+}$ be the subset of
dominant weights and $(\varpi_{i})_{i=1,\dots,r}$ be the fundamental
weights. Given a dominant weight $\lambda$, we denote
$L({\lambda})$ the associated simple $\g$-module with highest weight
$\lambda$. Finally
$(\,,\,)$ will be the Cartan-Killing form normalized such that for
the highest root
$\theta$ we have $(\theta,\theta)=2$. Let $\Lg=\g\otimes_{\comp}\comp((z))$ be the
{\it loop algebra} of 
$\g$ and $\Lgh$ be the central extension of $\Lg$
\begin{formula}\label{form:cent_ext}
0\lra\comp\lra\Lgh\lra\Lg\lra0
\end{formula}
defined by the $2$-cocycle 
$(X\otimes f,Y\otimes g)\mapsto (X,Y)\Res_{0}(gdf).$

\ind Fix an integer $\ell$. Call a representation of $\Lgh$ of level
$\ell$ if the center acts by multiplication by $\ell$. The theory of affine Lie
algebras affirms that the irreducible and integrable representations of $\Lgh$
are classified by the dominant weights belonging to $P_{\ell}=\{\lambda\in
P_{+}/(\lambda,\theta)\leq\ell\}$. For $\lambda\in P_{\ell}$, denote
${\cal{H}}_{\ell}(\lambda)$ the associated representation.

\subsection{\em Definition of conformal blocks. }\label{subsec:conf-blocks} 
Fix an integer (the level) $\ell\geq0$.
Let $(X,\ul{p})$ be an $n$-pointed stable curve (we denote 
$\ul{p}=(p_{1},\dots,p_{n})$) and suppose that the points are labeled by
$\ul{\lambda}=(\lambda_{1},\dots,\lambda_{n})\in P_{\ell}^{n}$
respectively. Choose a non-singular point $p\in X$ and a local coodinate $z$ at
$p$. Let
$X^{*}=X\moins\{p\}$ and $\LgX$ be the Lie algebra $\g\otimes{\cal{O}}(X^{*})$.  We
have a morphism on Lie algebras
$\LgX\ra\Lg$ by associating to $X\otimes f$ the element $X\otimes\hat{f}$, where 
$\hat{f}$ is the Laurent developpement of $f$ at $p$. By the
residue theorem, the restriction to $\LgX$ of the central extension
(\ref{form:cent_ext}) splits and we may see $\LgX$ as a Lie subalgebra of $\Lgh$. In
particuler, the
$\Lgh$-module ${\cal{H}}_{\ell}(0)$ may be seen as a $\LgX$-module. In addition, we
may consider the $\g$-modules $L(\lambda_{i})$ as a $\LgX$-modules by evaluation at
$p_{i}$. The vector space of conformal blocks is defined as
follows:
\begin{formula}\label{def:conf-blocks}
B_{G,X}(\ell;\ul{p};\ul\lambda)=[{\cal{H}}_{\ell}(0)\otimes_{\comp}
L(\lambda_{1})\otimes_{\comp}\dots_{\comp}L(\lambda_{n})]_{\LgX}
\end{formula}
where $[]_{\LgX}$ means that we take co-invariants. It is known 
(\cite{TUY:conf-field} or \cite{So:NB794}, 2.5.1) that these vector spaces are
finite-dimensional. Important properties are as follows: 

\ind $a)$ $\dim B_{G,\proj_{1}}(\ell;p_{1};0)=1$

\ind $b)$ If one
adds a non-singular point
$q\in X$, then the 
spaces $B_{G,X}(\ell;\ul{p};\ul\lambda)$ and
$B_{G,X}(\ell;\ul{p},q;\ul\lambda,0)$ are canonically isomorphic 
(\cite{So:NB794}, 2.3.2).

\ind $c)$ Suppose $X$ is singular in $c$ and let $\widetilde{X}\ra X$ be a
partial desingularization of $c$. Let  $a$ and $b$ be the points of $\widetilde{X}$
over $c$. Then there is a canonical isomorphism
$$\bigoplus_{\mu\in P_{\ell}}
B_{G,X}(\ell;\ul{p},a,b;\ul\lambda,\mu,\mu^{*})\iso
B_{G,X}(\ell;\ul{p};\ul\lambda)$$

\ind $d)$ The dimension of $B_{G,X}(\ell;\ul{p};\ul\lambda)$ does not change
when $(X;\ul{p})$ varies in the stack of $n$-pointed stable curves
${\goth{M}}_{g,n}$ (\cite{TUY:conf-field}).

\subsection{\em Application: }\label{subsec:E_8} Consider the case of $G=E_8$ and
level
$1$ and remark that  $P_{1}$ contains {\em only} the trivial representation. 
In order to calculate $B_{E_8,X}(\ell;p;0)$, one reduces to $\proj_{1}$ with
points labeled with the trivial representation using $c)$ and
$d)$, then it follows from $b)$ and $a)$ that it is one-dimensional.

\section{The Picard group of $\M$}

\subsection{}  We recall the description of $\Pic(\M)$ of
\cite{L-S:verlinde}, which uses as main tool the {\it
uniformization} theorem which I now recall. Let $\LG$ be the loop
group $G\bigl(\comp((z))\bigr)$, seen as an ind-scheme over $\comp$,
$\LGp$ the sub-group scheme $G\bigl(\comp[[z]]\bigr)$ and
$\Q=\LG/\LGp$ be the infinite Grassmannian, which is a direct limit
of projective integral varieties ($loc.\, cit.$). Finally let $\LGX$
be the sub-ind group $G({\cal{O}}(X^{*}))$ of
$\LG$. The uniformization theorem (\cite{L-S:verlinde}, 1.3)
states that there is a canonical isomorphism of stacks 
$\LGX\bk\Q\iso\M$ and moreover that $\Q\ra\M$ is a $\LGX$-bundle.

\ind Let $\Pic_{\LGX}(\Q)$ be the group of $\LGX$-linearized line
bundles on
$\Q$. Recall that a
$\LGX$-linearization of the line bundle ${\scr{L}}$ on $\Q$ is an
isomorphism
$m^{*}{\scr{L}}\iso pr_{2}^{*}{\scr{L}}$, where
$m:\LGX\times\Q\ra\Q$ is the action of $\LGX$ on $\Q$, satisfying the
usual cocycle condition. It follows from the uniformization theorem
that
$$\Pic(\M)\iso\Pic_{\LGX}(\Q),$$ hence in order to understand
$\Pic(\M)$ it suffices to understand
$\Pic_{\LGX}(\Q)$. The Picard group of $\Q$ itself is infinite
cyclic; let me recall how its positive generator may be defined in
terms of central extensions of $\LG$.

\subsection{} If ${\cal{H}}$ is an (infinite) dimensional vector
space over $\comp$, we define the $\comp$-space $\End({\cal{H}})$ by
$R\mapsto\End({\cal{H}}\otimes_{\comp}R)$, the $\comp$-group
$GL({\cal{H}})$ as the group of its units and 
$PGL({\cal{H}})$ by $GL({\cal{H}})/G_{m}$. The $\comp$-group
$\LG$ acts on $\Lg$ by the adjoint action which is extended to
$\Lgh$ by the following formula:
$$\Ad(\gamma).(\alpha^{\prime},s)=\bigl(\Ad(\gamma).\alpha^{\prime},
s+\Res_{z=0}(\gamma^{-1}\frac{d}{dz}\gamma,\alpha^{\prime})\bigr)$$
where $\gamma\in\LG(R)$, $\alpha=(\alpha^{\prime},s)\in\Lgh(R)$ and
$(\,,\,)$ is the $R((z))$-bilinear extension of the Cartan-Killing
form. The main tool we use is that if
$\bar\pi:\Lgh\ra\End({\cal{H}})$ is an integral highest weight
representation, then for $R$ a
$\comp$-algebra and
$\gamma\in\LG(R)$ there is, locally over $\Spec(R)$, an automorphism
$u_{\gamma}$ of ${\cal{H}}_{R}={\cal{H}}\otimes_{\comp}R$, unique up
to
$R^{*}$, such that
\begin{formula}\label{form:Faltings}
\begin{diagram} {\cal{H}}&\efl{\bar\pi(\alpha)}{}&{\cal{H}}\\
\sfl{u_{\gamma}}{}&&\sfl{}{u_{\gamma}}\\
{\cal{H}}&\efl{\bar\pi(\Ad(\gamma).\alpha)}{}&{\cal{H}}
\end{diagram}
\end{formula} is commutative for any $\alpha\in\Lgh(R)$
(\cite{L-S:verlinde}, Prop. 4.3).

\ind By the above, the representation $\bar\pi$ may be ``integrated"
to a (unique) {\it algebraic} projective representation of $\LG$,
\ie that there is a morphism of
$\comp$-groups
$\pi:\LG\ra PGL({\cal{H}})$ whose derivate coincides with $\bar\pi$
up to homothety. Indeed, thanks to the unicity property the
automorphisms $u$ associated locally to $\gamma$ glue together to
define an element
$\pi(\gamma)\in PGL({\cal{H}})(R)$ and still because of the unicity
property, $\pi$ defines a morphism of $\comp$-groups. The assertion
on the derivative is consequence of (\ref{form:Faltings}). We apply
this to the basic representation ${\cal{H}}_{1}(0)$ of $\Lgh$.
Consider the central extension
\begin{formula}\label{gl(H)-ext} 1\lra G_{m}\lra
GL({\cal{H}}_{1}(0))\lra PGL({\cal{H}}_{1}(0))\lra 1.
\end{formula} The pull back of (\ref{gl(H)-ext}) to
$\LG$ defines a central extension to which we refer as the {\it
canonical} central extension of $\LG$:
\begin{formula}\label{can-ext} 1\lra G_{m}\lra\LGh\lra\LG\lra 1
\end{formula}  A basic fact is that the extension (\ref{can-ext})
splits canonically over $\LGp$ (\cite{L-S:verlinde}, 4.9),
hence we may define a line bundle on the homogeneous space
$\Q=\widehat\LG/\widehat\LGp$ via the character
$G_m\times\LGp\ra G_m$ defined by the first projection. Then this
line bundle generates $\Pic(\Q)$ (\cite{L-S:verlinde}, 4.11);
we denote by
${\cal{O}}_{\Q}(1)$ its dual.

\subsection{} By (\cite{L-S:verlinde}, 6.2) the forgetful
morphism
$\Pic_{\LGX}(\Q)\ra\Pic(\Q)$ is injective, and moreover 
($loc.\,cit.$, 6.4), the line bundle ${\cal{O}}_{\Q}(1)$
admits a
$\LGX$-linearization if and only if the restriction  of the central
extension (\ref{can-ext}) to $\LGX$ splits.  It is shown in
\cite{L-S:verlinde} that this is indeed the case for classical $G$ and
$G_2$ by directly constructing line bundles on $\M$ which pull back
to
${\cal{O}}_{\Q}(1).$ In one case the existence of the splitting can
be proved directly:

\begin{th}{Proposition} The restriction  of the central extension
(\ref{can-ext}) to $\LGX$ splits for $G=E_8$.
\end{th} {\it Proof:} Let ${\cal{H}}={\cal{H}}_{1}(0)$. It suffices
to show that the representation 
$\bar\pi:\LgX\ra\End({\cal{H}})$ integrates to an algebraic
representation $\pi:\LGX\ra GL({\cal{H}})$, which in turn will follow
from the fact that in the case $\gamma\in\LGX(R)$ we can {\it
normalize} the automorphism $u_{\gamma}$ of (\ref{form:Faltings}).
Indeed, the commutativity of (\ref{form:Faltings}) shows that
coinvariants are mapped to coinvariants under
$u_{\gamma}$. For ${\goth{g}}=e_{8}$, $\ell=1$ and $\lambda=0$, we
know by (\ref{subsec:E_8}) that these spaces are
$1$-dimensional, hence we may choose $u_{\gamma}$ (in a unique way)
such that it induces the identity on coinvariants. \cqfd

\begin{th}{Corollary} Suppose $G=F_4,E_6,E_7$ or $E_8$.  There is a
line bundle ${\cal{L}}$ on $\M$ such that the pullback to $\Q$ is
${\cal{O}}_{\Q}(1)$.
\end{th} {Proof:} For $E_8$, this follows from the above
proposition.  Now consider the well known tower of natural inclusions
\begin{formula}\label{tower} F_{4}\rInto^{\alpha}_{}
E_6\rInto^{\beta}_{} E_7\rInto^{\gamma}_{} E_8.
\end{formula} On the level of Picard groups we deduce
$$\begin{diagram}[silent]
\Pic(\QEeight)&\rTo^{\tilde f_\alpha^{*}}_{}&\Pic(\QEseven)&
\rTo^{\tilde f_\beta^{*}}_{}&\Pic(\QEsix)&\rTo^{\tilde
f_\gamma^{*}}_{}&\Pic(\QFfour)\cr
\nfl{\pi_{E_8}^{*}}{}&&\nfl{\pi_{E_7}^{*}}{}&&\nfl{\pi_{E_6}^{*}}{}&&
\nfl{\pi_{F_4}^{*}}{}\cr
\Pic(\MEeight)&\rTo^{f_\alpha^{*}}_{}&\Pic(\MEseven)&
\rTo^{f_\beta^{*}}_{}&\Pic(\MEsix)&\rTo^{f_\gamma^{*}}_{}&\Pic(\MFfour)\cr
\nfl{f_{\varpi_{8}}^{*}}{}&&&&&&\ruTo(12.8,2.2)_{f_{\varpi_{8|F_4}}^{*}}\cr
\Pic(\MSLTF8)
\end{diagram}
$$ The Dynkin index of the representation $\varpi_{8}$ of $E_8$ is
60, and an easy calculation shows that 
$\varpi_{8|F_4}=14\,\comp\oplus \varpi_{1}\oplus 7\,\varpi_{4}$,
hence is equally of Dynkin index 60
(\cite{K-N:picard},\cite{L-S:verlinde}, 2.3). By the
Kumar-Narasimhan-Ramanathan lemma  (\cite{L-S:verlinde}, 6.8) the
determinant bundle 
${\cal{D}}$ pulls back, to
${\cal{O}}_{\QEeight}(60)$ via $\pi_{E_8}\circ f_{\varpi_{8}}^{*}$
and to ${\cal{O}}_{\QFfour}(60)$ via  $\pi_{F_4}\circ
f_{\varpi_{8|F_4}}^{*}$. If follows that $\tilde f_\alpha^{*},\tilde
f_\beta^{*}$ and $\tilde f_\gamma^{*}$ are isomorphisms and that the
pullback of the line bundle 
${\cal{L}}$ on $\MEeight$ under $f_\alpha^{*}$ (resp.
$f_\beta^{*}\circ f_\alpha^{*}$,
$f_\gamma^{*}\circ f_\beta^{*}\circ f_\alpha^{*}$) pulls back to
${\cal{O}}_{\QEseven}(1)$ (resp. ${\cal{O}}_{\QEsix}(1)$,
${\cal{O}}_{\QFfour}(1))$ \cqfd

\subsection{Proof of
theorem \ref{th:local_factoriality}: }
According to (\cite{B-L-S:picard}, 13)  it remains to prove that
$\Mod$ is not locally factorial for
$G=F_4,E_6,E_7$ or $E_8$. In order to see this we consider again the
tower (\ref{tower}) with additionally the natural inclusion
$\Spin_{8}\rInto^{}F_4$. Again the restriction of the
representation
$\varpi_{8}$ of
$E_8$ to
$\Spin_{8}$ has Dynkin index $60$, hence if the generator of
$\Pic(\M)$ would exist on $\Mod$, then
the Pfaffian bundle would exist on
$\ModSpinEight$, which is not the case (\cite{B-L-S:picard}, 8.2).
But the generators exist on the open subset of regularly stable
bundles, as the center of $G$ acts trivally on the fibers by
construction (we started with the trivial representation) and then the
arguments of (\cite{B-L-S:picard}, 13) apply. \cqfd

{\eightpoint
\bibliography{sorger}
\bibliographystyle{sorger} }

\bigskip
\hfill\vtop{\eightpoint
   \hbox to 5cm{\hfill Christoph {\sc Sorger}\hfill}
   \hbox to 5cm{\hfill Institut de math\'ematiques -- CP 7012\hfill} 
   \hbox to 5cm{\hfill Universit\'e Paris 7\hfill}
   \hbox to 5cm{\hfill 2, place Jussieu\hfill}
   \hbox to 5cm{\hfill F-75251 {\sc PARIS} Cedex 05\hfill}}

\end{document}